\documentclass[11pt]{article}

\usepackage{amsmath}
\usepackage{amsfonts}
\usepackage{amssymb}
\usepackage{amsthm}
\usepackage{graphics}
\usepackage{amscd}
\usepackage{graphicx,epsfig}
\usepackage{color}
\usepackage[all]{xy}
\usepackage{url}

\renewcommand{\epsilon}{\varepsilon}

\newcommand{\boL}{\mathcal{L}}

\newcommand{\boK}{\mathcal{K}}

\newcommand{\R}{\mathbb{R}}

\newcommand{\M}{\mathbb{M}}

\newcommand{\id}{\text{id}}

\newcommand{\dd}{\mathrm{d}}

\newtheorem{thm}{Theorem}
\newtheorem{cor}[thm]{Corollary}

\renewcommand{\phi}{\varphi}
\newcommand{\dis}{\displaystyle}

\newtheorem*{thm*}{Theorem}

\theoremstyle{remark}

\newtheorem*{rem*}{Remark}

\newcounter{remark}

\newcounter{case}

\newcounter{construction}

\newcounter{fact}

\DeclareMathOperator{\ric}{Ric}

\title{Optimal length estimates for stable CMC surfaces in $3$-space-forms}
\author{Laurent Mazet}
\date{}

\begin{document}

\maketitle

\begin{abstract}
In this paper, we study stable constant mean curvature $H$ surfaces in
$\R^3$. We prove that, in such a surface, the distance from a point to the
boundary is less that $\pi/(2H)$. This upper-bound is optimal and is
extended to stable constant mean curvature surfaces in space forms.
\end{abstract}

\section{Introduction}

A constant mean curvature (cmc) surface $\Sigma$ in a Riemannian
$3$-manifold
$\M^3$ is stable, if its stability operator $L=-\Delta-\ric(n,n)-|A|^2$ is
nonnegative. The nonnegativity of this operator means that $\Sigma$ is a
local minimizer of the area functional on surfaces regard to the
infinitesimal deformations fixing its boundary.

The stability hypothesis was studied by several authors and has many
consequences (see \cite{MePeRo} for an overview). For example,
D.~Fischer-Colbrie and R.~Schoen
\cite{FCSc} studied the case of complete stable minimal surfaces when
$\M^3$ has non-negative scalar curvature. They obtain that the universal
cover of $\Sigma$ is not conformally equivalent to the disk and, as a
consequence, prove that the plane is the only complete stable minimal
surface in $\R^3$. From this, R.~Schoen \cite{Sch} has derived a
curvature estimate for stable cmc surfaces.

In \cite{CoMi}, T.~H.~Colding and W.~P.~Minicozzi introduced new
technics and obtained area and curvature estimates for stable cmc surfaces.
Afterward, these technics were used by P.~Castillon \cite{Cas} to answer a
question asked in \cite{FCSc} about the consequences of the
positivity of certain elliptic operators. Recently, the same ideas have
been used by J.~Espinar and H.~Rosenberg \cite{EsRo} to obtain similar
results.

In \cite{RoRo}, A.~Ros and H.~Rosenberg study constant mean curvature $H$
surfaces in $\R^3$ with $H\neq 0$ : they prove a maximum principle at
infinity. One of their tools is a length estimate for stable cmc surface.
In fact, they prove that the intrinsic distance from a point $p$
in a stable cmc surface $\Sigma$ to the boundary of $\Sigma$ is less than
$\pi/H$. The aim of this paper is to improve this result. In fact,
applying the ideas of \cite{CoMi}, we prove that the distance is less than
$\pi/(2H)$. This estimate is optimal since, for a hemisphere of radius
$1/H$, the distance from the pole to the boundary is $\pi/(2H)$. Actually
we prove that the hemisphere of radius $1/H$ is the only stable cmc $H$
surface where the distance $\pi/(2H)$ is reached. We can generalized this
result to stable cmc
$H$ surfaces in $\M^3(\kappa)$, where $\M^3(\kappa)$ is the $3$-space form
of sectional curvature $\kappa$. We prove that when $H^2+\kappa>0$ such
an optimal estimate exists. In fact, it is already known that, when
$\kappa\le 0$ and $H^2+\kappa\le 0$, there is no such estimate since
 there exist complete stable cmc $H$ surfaces. But, in some sense, our
results is an extension of the fact that
the planes (resp. the horospheres) are the only stable complete constant
mean curvature $H$ surfaces in $\R^3$ (resp. $M^3(\kappa)$, $\kappa<0$)
when $H=0$ (resp. $H^2+\kappa=0$).

\section{Definitions}

On a constant mean curvature surface $\Sigma$ in a Riemannian $3$-manifold
$\M^3$, the stability operator is defined by $L=-\Delta
-\ric(n,n)-|A|^2$, where $\Delta$ is the Laplace operator
on $\Sigma$, $\ric$ is the Ricci tensor on $\M^3$, $n$ is the normal to
$\Sigma$ and $A$ is the second fundamental form on $\Sigma$.  When it is
necessary, we will denote the stability operator by $L_f$ to refer to the
immersion $f$ of $\Sigma$ in $\M^3$.

The surface $\Sigma$ is called \emph{stable} if the operator $L$ is
nonnegative \textit{i.e.}, for every compactly supported function $u$, we
have
$$
0\le \int_\Sigma u L(u)\dd \sigma=\int_\Sigma \|\nabla
u\|^2-(\ric(n,n)+|A|^2)u^2\dd
\sigma
$$
We remark that this property is sometimes called strong stablility since it
means that the second derivatives of the area functional is nonnegative
with
respect to any compactly supported infinitesimal deformations $u$ whereas
$\Sigma$ is critical for this functional only for compactly supported
infinitesimal deformations with vanishing mean value \textit{i.e.}
$\int_\Sigma u\dd \sigma=0$.

In the following, on a cmc surface, the normal $n$ is always chosen such
that $H$ is non-negative.

We will denote by $d_\Sigma$ the intrinsic distance on $\Sigma$ and by
$K$ the sectional curvature of the surface.

\section{Results}

The main result of this paper is the following theorem.

\begin{thm}\label{main}
Let $H$ be positive. Let $\Sigma$ be a stable constant mean curvature $H$
surface in $\R^3$.
Then, for $p\in\Sigma$, we have :
\begin{equation}\label{estimate}
d_\Sigma(p,\partial\Sigma)\le \frac{\pi}{2H}
\end{equation}
Moreover, if the equality is satisfied, $\Sigma$ is a hemisphere.
\end{thm}

In $\R^3$, the stability operator can be written $L=-\Delta-4H^2+2K$.

\begin{proof}
We denote by $R_0$ the distance $d_\Sigma(p,\partial\Sigma)$ and assume
that $R_0\ge\pi/(2H)$. If $R_0<\pi/H$ we denote by $I$ the segment
$[\pi/(2H), R_0]$, otherwise $I=[\pi/(2H),\pi/H)$. In fact, because of
the work of Ros and Rosenberg \cite{RoRo}, we already know that $R_0\le
\pi/H$. Let $R$ be in $I$.

The surface $\Sigma$ has constant mean curvature $H$ thus its sectional
curvature is less than $H^2$. So the exponential map $\exp_p$ is a local
diffeomorphism on the disk $D(0,R)\subset T_p\Sigma$ of center $0$ and
radius $R$. On this disk, we consider the induced metric and the operator
$\boL=-\Delta-4H^2+2K$. The surface $\Sigma$ is stable so it exists a
positive function $g$ on $\Sigma$ such that $L(g)=0$ (see
Theorem~1 in \cite{FCSc}). On $D(0,R)$, the
function $\tilde{g}=g\circ \exp_p$ is then positive and satisfies
$\boL(\tilde{g})=0$ since $D(0,R)$ and $\Sigma$ are locally isometric.
The operator $\boL$ is thus nonnegative on $D(0,R)$ \cite{FCSc}.

For $r\in[0,R]$, we define $l(r)$ as the length of the circle $\{v,\
 |v|=r\}\subset D(0,R)$ and $\boK(r)=\int_{D(0,r)}K\dd \sigma$. Since
$D(0,R)$
and $\Sigma$ are locally isometric, the sectional curvature $K$ of
$D(0,R)$ is less than $H^2$. Then
\begin{equation}\label{longueur}
l(r)\ge \frac{2\pi}{H}\sin Hr
\end{equation}
By Gauss-Bonnet, we have:
\begin{equation}\label{gaussbonnet}
\boK(r)=2\pi-l'(r)
\end{equation}

Let us consider a function $\eta:[0,R]\rightarrow[0,1]$ with $\eta(0)=1$
and $\eta(R)=0$. Let us write the nonnegativity of $\boL$ for the
radial function $u=\eta(r)$.
$$
0\le \int_0^R(\eta'(r))^2l(r)\dd r-4H^2\int_0^R\eta^2(r)l(r)\dd
r+2\int_0^R\boK'(r)\eta^2(r)\dd r
$$
Hence, following the ideas in \cite{CoMi} and using \eqref{gaussbonnet}
and the boundary values of $\eta$, we have:
\begin{align*}
\int_0^R(4H^2\eta^2-{\eta'}^2)l\dd r&\le 2\left(\left[\boK(r)\eta^2(r)
\right]_0^R -\int_0^R\boK(r)(\eta^2(r))'\dd r\right)\\
&\le -2\int_0^R\boK(r)(\eta^2(r))'\dd r\\
&\le -2\int_0^R(2\pi-l'(r))(\eta^2(r))'\dd r\\
&\le 4\pi+2\int_0^R(\eta^2(r))'l'(r)\dd r\\
&\le 4\pi+\left[2(\eta^2(r))'l(r)\right]_0^R-2\int_0^R(\eta^2(r))''l(r)
\dd r\\
&\le 4\pi-2\int_0^R(\eta^2(r))''l(r) \dd r
\end{align*}
Thus we obtain
\begin{equation}\label{positivity}
\int_0^R\left(4H^2\eta^2-{\eta'}^2+2(\eta^2)''\right)l\dd r\le 4\pi
\end{equation}
We shall apply this equation to the function $\eta(r)=\dis\cos\frac{\pi
r}{2R}$. In this case we have 
\begin{align*}
{\eta'}^2&=\frac{\pi^2}{4R^2}\sin^2\frac{\pi r}{2R}\\
(\eta^2)''&=-\frac{\pi^2}{2R^2}\left(\cos^2\frac{\pi
r}{2R}-\sin^2\frac{\pi r}{2R}\right)
\end{align*}
Thus 
$$
4H^2\eta^2-{\eta'}^2+2(\eta^2)''=(4H^2-\frac{\pi^2}{R^2})\cos^2\frac{\pi
r}{2R}+\frac{3\pi^2}{4R^2}\sin^2\frac{\pi r}{2R}
$$
As $R\ge \frac{\pi}{2H}$, $4H^2\eta^2-{\eta'}^2+2(\eta^2)''$ is
non-negative and, by \eqref{longueur},
{\footnotesize
\begin{align*}
\left(4H^2\eta^2-{\eta'}^2+2(\eta^2)''\right)l
&\ge\left((4H^2-\frac{\pi^2}{
R^2})\cos^2\frac{\pi r}{2R}+\frac{3\pi^2}{4R^2}\sin^2\frac{\pi
r}{2R}\right)\frac{2\pi}{H}\sin Hr\\
&\ge\frac{\pi}{H}\left((4H^2-\frac{\pi^2}{4R^2})\sin
Hr+(4H^2-\frac{7\pi^2}{4R^2})\frac{1}{2}\left(\sin(\frac{\pi}{R}+H)r-
\sin(\frac{\pi}{R}-H)r\right)\right)
\end{align*}
}

Thus integrating in \eqref{positivity}, we obtain (we recall that
$R<\pi/H$)
{\footnotesize
\begin{multline*}
4\pi\ge\frac{\pi}{H}\left((4H^2-\frac{\pi^2}{4R^2})\frac{1}{H}
(1-\cos HR)\right.\\
\left.+(4H^2-\frac{7\pi^2}{4R^2})\frac{1}{2}\left(\frac{R}{ \pi+HR }
(1-\cos(\pi+HR))-\frac{R}{\pi-HR}(1-\cos(\pi-HR))\right)\right)
\end{multline*}
}
After some simplifications in the above expression, we obtain
$$
4\pi\ge
\pi\frac{(-32H^2R^4+24\pi^2H^2R^2-\pi^4)-(10\pi^2H^2R^2-\pi^4)\cos
HR}{4H^2R^2(\pi^2-H^2R^2)}
$$
Now, passing $4\pi$ on the right-hand side of the above inequality and
simplifying by $\pi$, we obtain:
$$
0\ge \frac{-(4H^2R^2-\pi^2)^2-(10\pi^2H^2R^2-\pi^4)\cos
HR}{4H^2R^2(\pi^2-H^2R^2)}
$$
We denote by $F(R)$ the right-hand term of the above inequality. Hence we
have proved that, for every $R$ in $I$, $F(R)\le 0$. If we write
$R=\pi/(2H)+x$, we compute the Taylor
expansion of $F$ and obtain 
$$
F(\frac{\pi}{2H}+x)=2Hx+o(x)
$$
which is positive if $x>0$. Thus, if $R_0>\pi/(2H)$, we get a contradiction
and the inequality \eqref{estimate} is proved. 

Now if $R_0=\pi/(2H)$, we have in fact equality all along the computation,
so $l(r)=2\pi/H\sin Hr$ and $\boK(r)=2\pi-l'(r)=2\pi(1-\cos Hr)$. But we
also know
that the sectional curvature is less than $H^2$ thus $\boK(r)\le
H^2\int_0^r
l(u)\dd u=2\pi(1-\cos Hr)$. Since this inequality is in fact an equality,
the sectional curvature is in fact $H^2$ at every point. Thus the
principal curvatures of a point in $\Sigma$ are $H$ and $H$ \textit{i.e.}
there are only umbilical points. Hence $\Sigma$ is a piece of a sphere of
radius $1/H$ and, since $d_\Sigma(p,\partial\Sigma)= \frac{\pi}{2H}$, it
contains the hemisphere of pole $p$. A hemisphere can not be strictly
contained in a stable subdomain of the sphere, so $\Sigma$ is a
hemisphere.
\end{proof}

With this result we have an important corollary.

\begin{cor}Let $H\ge 0$ and $\kappa\in\R$ such that $H^2+\kappa>0$.
Let $\Sigma$ be a stable contant mean curvature $H$ surface in
$\M^3(\kappa)$. Then for $p\in\Sigma$, we have :
$$
d_\Sigma(p,\partial\Sigma)\le \frac{\pi}{2\sqrt{H^2+\kappa}}
$$
Moreover, if the equality is satisfied, $\Sigma$ is a geodesical
hemisphere of $\M^3(\kappa)$.
\end{cor}

The proof is based on the Lawson's correspondence between constant mean
curvature surfaces in space forms (see \cite{Law}).

\begin{proof}
First, the case $\kappa=0$ is Theorem~\ref{main}.

Let $\Pi:\widetilde{\Sigma}\rightarrow \Sigma$ be the universal cover of
$\Sigma$. We then have a constant mean curvature immersion of
$\widetilde{\Sigma}$ in $\M^3(\kappa)$, let $\boL=-\Delta -2\kappa-|A|^2$
be the stability operator on $\widetilde{\Sigma}$. $\Sigma$ is stable,
so there exists a positive function $g$ on $\Sigma$ such that
$L(g)=-\Delta g-(2\kappa+|A|^2)g=0$. Thus the function $\tilde{g}=g\circ
\Pi$ is a positive function on $\widetilde{\Sigma}$ satisfying
$\boL(\tilde{g})=0$. Hence $\widetilde{\Sigma}$ is stable. Let $\mathrm{I}$
and $S$
be respectively the first fundamental form and the shape operator
on $\widetilde{\Sigma}$. They satisfy the Gauss and Codazzi equations for
$\M^3(\kappa)$.

We define $S'=S+(-H+\sqrt{H^2+\kappa})\id$ on $\widetilde{\Sigma}$. Then
$\mathrm{I}$ and $S'$ satisfy the Gauss and Codazzi equations for
$\M^3(0)=\R^3$ (see \cite{Law}).
Hence there exists an immersion $f$ of $\widetilde{\Sigma}$ in $\R^3$ with
first fundamental form $\mathrm{I}$ and shape operator $S'$ (we notice that
the induced metric is the same). Its mean curvature is then
$H+(-H+\sqrt{H^2+\kappa})=\sqrt{H^2+\kappa}$ \textit{i.e.} the immersion
has constant mean curvature. The stability operator is
\begin{align*}
L_f&=-\Delta-\|S'\|^2\\
&=-\Delta-(\|S\|^2+4H(-H+\sqrt{H^2+\kappa } )+2(-H+\sqrt{H^2+\kappa})^2) \\
&=-\Delta-(\|S\|^2+2\kappa)\\
&=\boL
\end{align*}
Hence the surface $f(\widetilde{\Sigma})$ is stable. So, from
Theorem~\ref{main}, we have
$$
d_\Sigma(p,\partial\Sigma)=
d_{\widetilde{\Sigma}}(\tilde{p},\partial\widetilde{\Sigma})\le
\frac{\pi}{2\sqrt{H^2+\kappa}}
$$
where $\Pi(\tilde{p})=p$. 

The equality case comes from the equality case in Theorem~\ref{main} and
since the Lawson's correspondence sends spheres into spheres.
\end{proof}

\bibliographystyle{plain}
\bibliography{../reference}

%
%
%
%
%

\end{document}